# Proof of the gradient conjecture of R. Thom

By Krzysztof Kurdyka, Tadeusz Mostowski, and Adam Parusiński

**Abstract**

Let $x(t)$ be a trajectory of the gradient of a real analytic function and suppose that $x_0$ is a limit point of $x(t)$. We prove the gradient conjecture of R. Thom which states that the secants of $x(t)$ at $x_0$ have a limit. Actually we show a stronger statement: the radial projection of $x(t)$ from $x_0$ onto the unit sphere has finite length.

## 0. Introduction

Let $f$ be a real analytic function on an open set $U \subset \mathbb{R}^n$ and let $\nabla f$ be its gradient in the Euclidean metric. We shall study the trajectories of $\nabla f$, i.e. the maximal curves $x(t)$ satisfying

$$\frac{dx}{dt}(t) = \nabla f(x(t)),\ t \in [0, \beta).$$

In the sixties Łojasiewicz [Lo2] (see also [Lo4]) proved the following result.

Łojasiewicz's Theorem. *If $x(t)$ has a limit point $x_0 \in U$, i.e. $x(t_\nu) \to x_0$ for some sequence $t_\nu \to \beta$, then the length of $x(t)$ is finite; moreover, $\beta = \infty$. Therefore $x(t) \to x_0$ as $t \to \infty$.*

Note that $\nabla f(x_0) = 0$, since otherwise we could extend $x(t)$ through $x_0$. The purpose of this paper is to prove the following statement, called "the gradient conjecture of R. Thom" (see [Th], [Ar], [Lo3]):

Gradient Conjecture. *Suppose that $x(t) \to x_0$. Then $x(t)$ has a tangent at $x_0$, that is the limit of secants $\lim\limits_{t\to\infty} \dfrac{x(t) - x_0}{|x(t) - x_0|}$ exists.*

This conjecture can be restated as follows. Let $\tilde{x}(t)$ be the image of $x(t)$ under the radial projection $\mathbb{R}^n \setminus \{x_0\} \ni x \longrightarrow \dfrac{x - x_0}{|x - x_0|} \in S^{n-1}$. The conjecture claims that $\tilde{x}(t)$ has a limit.

Actually we shall prove a stronger result: the length of $\tilde{x}(t)$ can be uniformly bounded for trajectories starting sufficiently close to $x_0$. In the last chapter we prove that the conjecture holds in the Riemannian case. More precisely, let $\nabla_g f$ be the gradient of $f$ with respect to some analytic Riemannian metric $g$ on $U$, suppose that $x(t)$ is a trajectory of $\nabla_g f$ and $x(t) \to x_0$; then $x(t)$ has a tangent at $x_0$.

The paper is organized as follows:

In Section 1 we recall the main argument in Łojasiewicz's theorem, we derive from it a notion of control function and we explain the crucial role it plays in the proof of the conjecture.

In Section 2 we recall known results and basic ideas about the conjecture; in particular we sketch the main idea of [KM], where the first proof of the conjecture was given. This section contains also some heuristic arguments which help to explain the construction of a control function. In the end we state some stronger conjectures on the behavior of a trajectory $x(t)$ near its limit point $x_0$.

Section 3 contains a detailed plan of the proof of the conjecture, i.e., the construction of a control function.

The proof of the conjecture is given in Sections 4–7.

In Section 8 we show that actually there is a uniform bound: the radial projections of all trajectories, having 0 as the limit point, have their lengths bounded by a universal constant.

Let $\tilde{x}_0$ denote the limit point of $\tilde{x}(t)$. In the second part of Section 8 we compare the distance $|\tilde{x}(t) - \tilde{x}_0|$ and $r = |x(t) - x_0|$. By the conjecture, $|\tilde{x}(t) - \tilde{x}_0| \to 0$ as $r \to 0$, but, as some examples show, it may go to 0 more slowly than any positive power of $r$. Geometrically this means that there is no "cuspidal neighborhood" of the tangent line $P$ to $x(t)$ at $x_0$, of the form $\{x \in \mathbb{R}^n; \operatorname{dist}(x, P) \leq r^\delta\}$, $\delta > 1$, which captures the trajectory near the limit point.

Finally in Section 9 we prove the gradient conjecture in the Riemannian case by reducing it to the Euclidean case.

*Notation and conventions.* In the sequel we shall always assume $x_0 = 0$ and $f(0) = 0$, so that in particular, $f$ is negative on $x(t)$. We often write $r$ instead of $|x|$ which is the Euclidean norm of $x$. We use the standard notation $\varphi = o(\psi)$ or $\varphi = O(\psi)$ to compare the asymptotic behavior of $\varphi$ and $\psi$, usually when we approach the origin. We write $\varphi \sim \psi$ if $\varphi = O(\psi)$ and $\psi = O(\varphi)$, and $\varphi \simeq \psi$ if $\frac{\varphi}{\psi}$ tends to 1.

## 1. Łojasiewicz's argument and control functions

We shall usually parametrize $x(t)$ by its arc-length $s$, starting from point $p_0 = x(0)$, and

$$\dot{x} = \frac{dx}{ds} = \frac{\nabla f}{|\nabla f|}.$$

By Łojasiewicz's theorem the length of $x(s)$ is finite. Denote it by $s_0$. Then $x(s) \to 0$ as $s \to s_0$.

Our proof is modeled on Łojasiewicz's idea [Lo2] so we recall first his argument. The key point of this argument is the Łojasiewicz inequality for the gradient [Lo1] which states that in a neighborhood $U_0$ of the origin

(1.1) $$|\nabla f| \geq c|f|^\rho$$

for some $\rho < 1$ and $c > 0$. Thus in $U_0$ we have on the trajectory $x(s)$

(1.2) $$\frac{df}{ds} = \langle \nabla f, \dot{x} \rangle = |\nabla f| \geq c|f|^\rho.$$

In particular $f(x(s))$ is increasing and

$$\frac{d|f|^{1-\rho}}{ds} \leq -[c(1-\rho)] < 0,$$

the sign coming from the fact that $|f|$ is decreasing on the trajectory. The integration of $|f|^{1-\rho}$ yields the following: if $x(s)$ lies in $U_0$ for $s \in [s_1, s_2]$, then the length of the segment of the curve between $s_1$ and $s_2$ is bounded by

$$c_1[|f(x(s_1))|^{1-\rho} - |f(x(s_2))|^{1-\rho}],$$

where $c_1 = [c(1-\rho)]^{-1}$. Consequently, if the starting point $p_0 = x(0)$ is sufficiently close to the origin, then:

1. The length of $x(s)$ between $p_0 = x(0)$ and the origin is bounded by $c_1|f(p_0)|^{1-\rho}$.

2. The curve cannot leave $U_0$.

3. $|f(x(s))| \geq c_2 r^N$, where $r = |x(s)|$, $N = 1/(1-\rho)$, $c_2 = c_1^{-N}$.

In this paper we shall often refer to the argument presented above as *Łojasiewicz's argument*.

A *control function*, say $g$, for a trajectory $x(s)$, is a function defined on a set which contains the trajectory, such that $g(x(s))$ grows "fast enough." In Łojasiewicz's proof, the function $f$ itself is a control function; what "grows fast enough" means is given by (1.2); it is this rate of growth, together with boundedness of $f$, which implies that the length of $x(s)$ is finite.

To illustrate how a control function will be used, consider the radial projection $\tilde{x}(s)$ of the trajectory. Let us parametrize $\tilde{x}(s)$ by its arc-length $\tilde{s}$. We use $\tilde{s}$ to parametrize the trajectory itself. Assume that we have a control function $g$, bounded on the trajectory, such that the function $\tilde{s} \to g(x(\tilde{s}))$ grows sufficiently fast. Then the length of $\tilde{x}(s)$ must be finite, as in Łojasiewicz's argument. As we show in Section 7, for a given trajectory $x(t)$ there exists a control function of the form

$$g = F - a - r^\alpha,$$

where $a$ is a negative constant, $\alpha > 0$ is small enough and $F = \frac{f}{r^l}$ with some rational $l$.

More precisely, $g$ is bounded on the trajectory and satisfies $\frac{dg}{d\tilde{s}} \geq |g|^\xi$, for a $\xi < 1$. Hence the proof can be completed by the Łojasiewicz argument.

## 2. Geometric motivations and historical account

We shall discuss now some known cases of the gradient conjecture and some ideas related to its proof.

Let us expand the function $f$ in the polar coordinates $(r, \theta)$ in $\mathbb{R}^n$, with $\theta \in S^{n-1}$:

(2.1) $$f = f_0(r) + r^m F_0(\theta) + \ldots,$$

$F_0 \neq \text{const}$. If $m = \infty$, then all trajectories of $\nabla f$ are straight lines. Now the equations $\frac{dx}{dt}(t) = \nabla f$ in polar coordinates are

(2.2) $$\frac{dr}{dt} = \frac{\partial f}{\partial r}, \qquad \frac{d\theta}{dt} = r^{-2} \frac{\partial f}{\partial \theta}.$$

The spherical part of $f$, i.e. $F_0(\theta)$, can be considered as a function on $S^{n-1}$ or as a function on $\mathbb{R}^n \setminus \{0\} \ni x \mapsto F_0\left(\frac{x}{|x|}\right)$.

If the order $d$ of $f_0(r)$ is smaller than or equal to $m - 1$, then the gradient conjecture is easy since for some $C > 0$

$$\left|\frac{d\theta(x(r))}{dr}\right| < C r^{m-d-1} < C,$$

for $r = |x(s)| < 1$. Now the length of $\tilde{x}(r) = \theta(x(r))$ is finite since the length of $r(x(s)) = |x(s)|$ is finite.

R. Thom, J. Martinet and N. Kuiper (see also F. Ichikawa [Ic]) proved two cases of the gradient conjecture by using (2.2) and applying Łojasiewicz's argument to $F_0$. They proved that:

1. $F_0(x(s))$ has a limit $\alpha \leq 0$,

2. if $\alpha < 0$ or if $\{F_0 = 0\} \cap \{\nabla F_0 = 0\}$ has only isolated points, then the limit of secants of $x(s)$ exists.

The proofs are published in [Mu].

*Attempts to construct a control function.* Let us again use expansion (2.1). Denote by $\nabla F_0$ the gradient $F_0$ with respect to the Riemannian metric induced on the sphere $S^{n-1}$ by the Euclidean metric on $\mathbb{R}^n$.

Assume first that $f$ is a homogeneous polynomial of degree $m$, that is $f = r^m F_0(\theta)$, $F_0 \neq 0$. It was observed by R. Thom that $\tilde{x}(s)$ is a trajectory of $\frac{\nabla F_0}{|\nabla F_0|}$; hence the gradient conjecture holds in this case. Moreover, by the Lojasiewicz inequality (1.1) for the function $F_0 : S^{n-1} \to \mathbb{R}$ it can be easily seen that $F_0 = \frac{f}{r^m}$ as a function on $\mathbb{R}^n \setminus 0$ is a control function for $x(\tilde{s})$.

In the general case it is easy to start to construct a control function. By the above result of Thom and Martinet we may assume that $\alpha = \lim_{s \to s_0} F_0(x(s)) = 0$. Suppose that $C > 0$ is big enough; then $F_0$ increases on $x(s)$ outside the set
$$\Omega_0 = \{x = (r, \theta) : |\nabla F_0(\theta)| < Cr\}.$$
Thus, $F_0$ may be considered as a control function, but only in $\mathbb{R}^n \setminus \Omega_0$. Using this fact we see that $x(s)$ must fall into $\Omega_0$ and cannot leave it.

Now it is natural to replace $\mathbb{R}^n$ by $\Omega_0$ and to try to construct a control function in $\Omega_0 \setminus \Omega_1$, where $\Omega_1$ is a (proper) subset of $\Omega_0$, and to prove that $x(s)$ must fall into $\Omega_1$, etc. More precisely we want to obtain a sequence $\Omega_i \subset \Omega_{i+1}$ such that dimensions of tangent cones (at 0) are decreasing.

Already the second step is not easy. Attempts to realize it were undertaken by N. Kuiper and Hu Xing Lin in the 3-dimensional case. Under an additional assumption Hu [Hu] succeeded in proving the gradient conjecture along these lines.

The first proof of the gradient conjecture in the general case was given in [KM]. Its starting point was that one can guess the end of the story of $\Omega_0, \Omega_1, \ldots$ . Indeed, let $x(s)$ be a trajectory of $\nabla f$ tending to 0.

For any $\lambda \in \mathbb{R}$ we define
$$(2.3) \qquad \mathcal{D}(\lambda) = \{x : |\nabla f(x)| < r^\lambda\}$$
and we put $k = \sup\{\lambda : \gamma \text{ intersects } \mathcal{D}(\lambda) \text{ in any nbd. of } 0\}$. We fix a rational $\lambda < k$ sufficiently close to $k$ and consider $\mathcal{D}(\lambda)$ as the "last" of $\Omega_i$.

A rather detailed analysis of the structure of $\mathcal{D}(\lambda)$ was done in the spirit of L-regular decomposition into L-regular sets (called also pancakes) of [Pa1], [Ku1], [Pa2]. As a result it was found that a function of the form
$$(2.4) \qquad F = \frac{f + cr^{k+1+\delta}}{r^{k+1}}$$

(for suitable constants $c, \delta$) should be taken as a control function. Its behavior was studied by different means inside $\mathcal{D}(\lambda)$ and outside $\mathcal{D}(\lambda)$. We didn't succeed in proving that $F$ increases always on $x(s)$, but we proved that it increases fast enough in most parts of $x(s)$. The final result was that the limit of secants exists.

*Blowing-up and the finiteness conjecture.* As far as we know the most common method (suggested also by R. Thom [Th]) to solve the gradient conjecture was to blow-up. Consider $p : M \to \mathbb{R}^n$ the blowing up of 0 in $\mathbb{R}^n$. Let $x^*(s)$ be the lifting of $x(s)$ via $p$; what one needs to prove the conjecture is that $x^*(s)$ has a limit as $s \to s_0$. One may try to follow Łojasiewicz; $x^*(s)$ is a trajectory of the gradient of $f \circ p$ in the "metric" induced by $p$; however, this "metric" degenerates on $p^{-1}(0)$ and the Łojasiewicz inequality (1.1) does not hold.

One may generalize this approach as follows. Let $V$ be any subvariety of the singular locus of $f$ (i.e. $df = 0$ on $V$), let $M$ be an analytic manifold and $p : M \to \mathbb{R}^n$ a proper analytic map such that $p : M \setminus p^{-1}(V) \to \mathbb{R}^n \setminus V$ is a diffeomorphism. For example, $p$ may be a finite composition of blow-ups with smooth centers. One may conjecture that the lifting of $x(t)$ to $M$ has a limit. Actually this follows from a stronger statement called the finiteness conjecture (proposed by R. Moussu and the first named author independently):

THE FINITENESS CONJECTURE FOR THE GRADIENT. *Let $A$ be a subanalytic subset of $\mathbb{R}^n$, then the set $\{t \in [0, \infty); x(t) \in A\}$ has finitely many connected components.*

Actually we can also consider an apparently weaker conjecture with a smaller class of sets, assuming that $A$ is an analytic subset of $\mathbb{R}^n$.

THE ANALYTIC FINITENESS CONJECTURE FOR THE GRADIENT. *Let $A$ be an analytic subset of $\mathbb{R}^n$; then either $x(t)$ stays in $A$ or it intersects $A$ in a finite number of points.*

The analytic finiteness conjecture implies (cf. [Ku2]) that the limit of tangents $\lim_{s \to s_0} \frac{x'(s)}{|x'(s)|}$ exists, which is still an open question in general, and which implies the gradient conjecture. Another consequence of the analytic finiteness conjecture is the positive answer to a conjecture of R. Thom that $|x(t)|$ is strictly decreasing from a certain moment.

F. Sanz [Sa] proved the analytic finiteness conjecture for $n = 3$ under the assumption that corank $D^2 f(0) = 2$. At the end of this section we propose a simple proof of the finiteness conjecture for $n = 2$. Recall that in this case any subanalytic set is actually semianalytic. Now, both finiteness conjectures are equivalent.

PROPOSITION 2.1. *If $\Gamma$ is an analytic subset in a neighborhood of $0 \in \mathbb{R}^2$ and $x(t) \to 0$ is a trajectory of $\nabla f$, then either $x(t)$ lies entirely in $\Gamma$ or it intersects $\Gamma$ in a finite number of points.*

*Proof.* First we show that $x(t)$ cannot spiral around the origin. For this purpose we expand the function $f$ in the polar coordinates $(r, \theta)$ in $\mathbb{R}^2$, so that we have as in (2.1)
$$f = f_0(r) + r^m F_0(\theta) + \ldots,$$
where $F_0 \neq \mathrm{const}$ is a function of $\theta \in \mathbb{R}$. If $m = \infty$, then all trajectories of $\nabla f$ are straight lines. Assume that $m < \infty$, then for some $\varepsilon > 0$ one of the sectors
$$A_+(\varepsilon) = \{x = (r, \theta) : F_0'(\theta) > \varepsilon\}, \qquad A_-(\varepsilon) = \{x = (r, \theta) : F_0'(\theta) < -\varepsilon\}$$
is not empty, and therefore, by periodicity of $F_0$, so is the other one (maybe for a smaller $\varepsilon$). It follows by (2.2) that, in a sufficiently small neighborhood of 0, the trajectory $x(t)$ crosses $A_+(\varepsilon)$ only anti-clockwise and $A_-(\varepsilon)$ clockwise. So $\theta$ is bounded on $x(t)$ and, in other words, the trajectory cannot spiral.

To end the proof of proposition take any semianalytic arc $\Gamma_1 \subset \Gamma \setminus \{0\}$, $0 \in \overline{\Gamma}_1$. If $\nabla f$ is tangent to $\Gamma_1$ and $x(t)$ meets $\Gamma_1$, then of course $x(t)$ stays in $\Gamma_1$. In the other case, $\nabla f$ is nowhere tangent to $\Gamma_1$, in a small neighborhood of 0. So $\Gamma_1$ can be crossed by $x(t)$ only in one way. Since $\theta$ is bounded on $x(t)$ and $\Gamma_1$ has tangent at 0, the trajectory meets $\Gamma_1$ only in a finite number of points. $\square$

## 3. The plan of the proof

The present proof is a simplified and modified version of the proof in [KM] proposed by the third named author. We shall outline below its main points. The proof is fairly elementary and is based on the theory of singularities.

First we replace the sets $D(\lambda)$, see (2.3), by much simpler sets $W^\varepsilon = \{x; f(x) \neq 0, \varepsilon|\nabla' f| \leq |\partial_r f|\}$, $\varepsilon > 0$, and then guess what the exponents $l = k+1$ of the denominators of (2.4) are. The role of $W^\varepsilon$ can be explained as follows. Decompose the gradient $\nabla f$ into its radial $\partial_r f \frac{\partial}{\partial r}$ and spherical $\nabla' f$ components, $\nabla f = \partial_r f \frac{\partial}{\partial r} + \nabla' f$, $r$ stands for $|x|$. Then, most of the time along the trajectory, the radial part must dominate; otherwise the trajectory would spiral and never reach the origin. Thus the trajectory $x(s)$ cannot stay away from the sets $W^\varepsilon$. On the contrary it has to pass through $W^\varepsilon$ in any neighborhood of the origin; for $\varepsilon > 0$ sufficiently small, see Proposition 6.2 below. But the limits of $\frac{r\partial_r f}{f}(x)$, as $x \to 0$, and $x \in W^\varepsilon$, are rational numbers and so form a finite subset $L$ of $\mathbb{Q}$; see Proposition 4.2. We call $L$ the set of characteristic exponents of $f$. It can be understood as a generalization

of the Łojasiewicz exponent; see Remark 4.4 below. Then, as we prove in Proposition 6.2, for each trajectory $x(s)$ there is an $l \in L$ such that $\frac{|f|}{r^l}(x(s))$ stays bounded from 0 and $\infty$. Thus $F = \frac{f}{r^l}$ is a natural candidate for a control function.

Let us remark that the main difficulty in proving the gradient conjecture comes from the movement of $x(s)$ in the sets of the form $\{r^\delta |\nabla' f| \leq |\partial_r f| \leq r^{-\eta}|\nabla' f|\}$, for $\delta > 0$ and $\eta > 0$. If $|\partial_r f| < r^\delta |\nabla' f|$, then the spherical part of the movement is dominant. Therefore not only $F$ but even $f$ itself can be used as a control function to bound the length of $\tilde{x}(s)$ (see, for instance, the last part of the proof of Theorem 8.1). On the other hand, if $r^{-\eta}|\nabla' f| < |\partial_r f|$, then the movement in the radial direction dominates. The function $-r^\alpha$, for any $\alpha > 0$, can be chosen as a control function; see the proof of Theorem 7.1.

Let us come back to $F = \frac{f}{r^l}$. We observe that $F(x(s))$ has a limit $a < 0$ as $x(s) \to 0$. This follows from the theory of asymptotic critical values of $F$ which we recall briefly in Section 5. Moreover, this limit has to be an asymptotic critical value of $F$ and the set of such values is finite. Finally, we have a strong version of the Łojasiewicz inequality for $F$, $r|\nabla F| \geq |F|^\rho$, $0 < \rho < 1$, not everywhere but at least on the sets where the main difficulty arises, that is, on the set $\{|\partial_r f| \leq r^{-\eta}|\nabla' f|\}$; see Proposition 5.3 and Lemma 7.2 below. Now the proof of the gradient conjecture is fairly easy. As we show in Section 7, $g = F - a - r^\alpha$, for $\alpha > 0$ and sufficiently small, is a good control function: it is bounded on the trajectory and satisfies $\frac{dg}{ds} \geq |g|^\xi$, for a $\xi < 1$. Hence the proof can be completed by the Łojasiewicz argument.

The main tools of the proof are the curve selection lemma and the classical Łojasiewicz inequalities. Only the proofs of Propositions 5.1 and 5.3 use the existence of Whitney stratification (with the (b) or (w) condition) of real analytic sets. For the existence of Whitney stratification see for instance [Lo1], [V]. A short and relatively elementary proof was presented also in [LSW]. Unlike the proof in [KM] our proof does not use L-regular sets, though it would be right to say that the study of L-regular decompositions led us, to a great extent, to the proof presented in this paper.

3.1. *A short guide on constants and exponents.* There are many equations and inequalities in the proof and each of them contains exponents and constants. Let us explain briefly the role of the most important ones. The constants are not important in general except for two of them: $c_f$, $\varepsilon$. The other constants just exist and usually we denote by $c$ the positive constants that are supposed to be sufficiently small and by $C$ the ones which are supposed to be sufficiently big. By $c_f$ we denote the constant of the Bochnak-Łojasiewicz inequality; see Lemma 4.3 below. For $c_f$ we may take any number smaller than the multiplicity of $f$ at the origin. For $\varepsilon$ we may take any positive number smaller than $\frac{1}{2}c_f(1 - \rho_f)$, where $\rho_f$ is the Łojasiewcz exponent, the smallest number $\rho$ satisfying (1.1). The set $L$ of characteristic exponents of $f$ is defined in the following section. For each exponent $l \in L$ there is an $\omega > 0$ defined in (6.4) (related to $\delta$ of Proposition 4.2). In general the letter $\delta$ may signify different exponents at different places of the proof, similarly to $c$ and $C$, but the exponent $\omega$ satisfying (6.4) is fixed (common for all $l \in L$ for simplicity). The other important exponent is $\alpha < \omega$, $\alpha > 0$, which is used in the formula (7.3) for the control function $g$. The exponent $\eta$ of the proof of the main theorem has auxiliary meaning, it allows us to decompose the set $W_l^\varepsilon$ into two pieces; $W_{-\eta,l}$ and to $W_l^\varepsilon \setminus W_{-\eta,l}$ and use different arguments on each piece; $\eta$ is chosen so that $\alpha < \eta < \omega$.

## 4. Characteristic exponents

Fix an exponent $\rho < 1$ so that in a neighborhood of the origin we have the Łojasiewicz inequality (1.1). The gradient $\nabla f$ of $f$ splits into its radial component $\frac{\partial f}{\partial r}\frac{\partial}{\partial r}$ and the spherical one $\nabla' f = \nabla f - \frac{\partial f}{\partial r}\frac{\partial}{\partial r}$. We shall denote $\frac{\partial f}{\partial r}$ by $\partial_r f$ for convenience.

For $\varepsilon > 0$ we denote $W^\varepsilon = \{x; f(x) \neq 0, \varepsilon|\nabla' f| \leq |\partial_r f|\}$. Note that $W^\varepsilon \subset W^{\varepsilon'}$ for $\varepsilon' < \varepsilon$.

LEMMA 4.1. *For each $\varepsilon > 0$, there exists $c > 0$, such that*

$$|f| \geq cr^{(1-\rho)^{-1}}, \tag{4.1}$$

*on $W^\varepsilon$. In particular each $W^\varepsilon$ is closed in the complement of the origin.*

*Proof.* Fix $\varepsilon > 0$. By the curve selection lemma it suffices to show that $|f|r^{-(1-\rho)^{-1}}$ is bounded from zero along any real analytic curve $\gamma(t) \to 0$ as $t \to 0$, $\gamma(t) \in W^\varepsilon$ for $t \neq 0$. Fix such a $\gamma$. In order to simplify the notation we reparametrize $\gamma$ by the distance to the origin, that is to say, $|\gamma(t(r))| = r$. In the spherical coordinates we write $\gamma(r) = r\theta(r)$, $|\theta(r)| = 1$. Then the tangent vector to $\gamma$ decomposes as the sum of its radial and spherical components as follows:

$$\gamma'(r) = \theta(r) + r\theta'(r),$$

and $r\theta'(r) = o(1)$. We have a Puiseux expansion

$$f(\gamma(r)) = a_l r^l + \ldots, \quad a_l \neq 0, \tag{4.2}$$

where $l \in \mathbb{Q}^+$, and

$$\frac{df}{dr}(\gamma(r)) = \partial_r f + \langle \nabla' f, r\theta'(r) \rangle = la_l r^{l-1} + \ldots. \tag{4.3}$$

By the assumption, $|\partial_r f| \geq \varepsilon |\nabla' f|$ on $\gamma$ and hence $\partial_r f$ is dominant in the middle term of (4.3). Consequently, along any real analytic curve $\gamma(r)$ in $W^\varepsilon$,

(4.4) $$r|\nabla f| \sim r|\partial_r f| \sim r|df/dr| \sim r^l \sim |f|.$$

In particular, by (1.1), $\frac{l-1}{l} \leq \rho$. This is equivalent to $(1-\rho)^{-1} \geq l$, and hence $|f|r^{-(1-\rho)^{-1}}$ is bounded from zero on $\gamma$ as required. □

Suppose we want to study the set of all possible limits of $\frac{r\partial_r f}{f}(x)$, as $W^\varepsilon \ni x \to 0$. By the curve selection lemma it suffices to consider real analytic curves $\gamma(r) \to 0$ contained in $W^\varepsilon$. For such a curve $\gamma(r)$, by (4.2) and (4.3), $\frac{r\partial_r f}{f} \to l$, where $l$ is a positive rational defined by (4.2). As we show below the sets of such possible limits is a finite subset $L \subset \mathbb{Q}^+$. By abuse of notation we shall write this property as $\frac{r\partial_r f}{f}(x) \to L$ for $W^\varepsilon \ni x \to 0$.

PROPOSITION 4.2.  *There exists a finite subset of positive rationals $L = \{l_1, \ldots, l_k\} \subset \mathbb{Q}^+$ such that for any $\varepsilon > 0$*

$$\frac{r\partial_r f}{f}(x) \to L \quad as \quad W^\varepsilon \ni x \to 0.$$

*In particular, as a germ at the origin, each $W^\varepsilon$ is the disjoint union*

$$W^\varepsilon = \bigcup_{l_i \in L} W^\varepsilon_{l_i},$$

*where we may define $W^\varepsilon_{l_i} = \{x \in W^\varepsilon; |\frac{r\partial_r f}{f} - l_i| \leq r^\delta\}$, for $\delta > 0$ sufficiently small.*

*Moreover, there exist constants $0 < c_\varepsilon < C_\varepsilon$, which depend on $\varepsilon$, such that*

(4.5) $$c_\varepsilon < \frac{|f|}{r^{l_i}} < C_\varepsilon \quad on \quad W^\varepsilon_{l_i}.$$

*Proof.* First we show that the set of possible limits is finite and independent of $\varepsilon$. Roughly speaking, the argument is the following. The set of limits of $\frac{r\partial_r f}{f}$, as $r \to 0$, is subanalytic, and hence, if contained in $\mathbb{Q}$, finite. We denote it by $L_\varepsilon$. Clearly $L_\varepsilon \subset L_{\varepsilon'}$ for $\varepsilon \geq \varepsilon'$. Moreover since $L_\varepsilon$, $\varepsilon \in \mathbb{R}^+$, is a subanalytic family of finite subanalytic subsets it has to stabilize, that is, $L_\varepsilon = L_{\varepsilon'}$ for some $\varepsilon > 0$ and each $0 < \varepsilon' \leq \varepsilon$.

We shall present this argument in more detail. Letting $\Omega = \{(x, \varepsilon); x \in W^\varepsilon, \varepsilon > 0\}$, we consider the map $\psi : \Omega \to \mathbb{P}^1$, where $\mathbb{P}^1 = \mathbb{R} \cup \{\infty\}$, defined by

$$\psi(x, \varepsilon) = \frac{r\partial_r f}{f}(x).$$

For any fixed $\varepsilon > 0$ the set $L_\varepsilon$, of limits of $\psi(x, \varepsilon)$ as $W^\varepsilon \ni x \to 0$, is a subanalytic set in $\mathbb{P}^1$ (in different terminology $L_\varepsilon \subset \mathbb{R}$ is subanalytic at infinity).

But by (4.4) and the curve selection lemma $L_\varepsilon \subset \mathbb{Q}^+$, hence $L_\varepsilon$ is finite. By a standard argument of subanalytic geometry the set

$$P = \{(\varepsilon, l); l \in L_\varepsilon, \varepsilon > 0\}$$

is subanalytic in $\mathbb{P}^1 \times \mathbb{P}^1$. Indeed, $P$ is obtained by taking limits at $0 \in \mathbb{R}^n$, with respect to $x$ variable, of a subanalytic function $\psi$. Since every $L_\varepsilon$ is finite there exists a finite partition $0 = \varepsilon_0 < \varepsilon_1 < \ldots < \varepsilon_N = +\infty$ such that $P \cap ((\varepsilon_i, \varepsilon_{i+1}) \times \mathbb{R})$ is a finite union of graphs of continuous functions on $(\varepsilon_i, \varepsilon_{i+1})$. But these functions take only rational (and positive) values; hence they are constant on $(\varepsilon_i, \varepsilon_{i+1})$. So $P \subset \mathbb{R} \times L$ for some finite subset $L$ of $\mathbb{Q}^+$, and we take $L$ to be the smallest with this property.

*Remark.* Actually in the sequel we shall work with $\varepsilon \to 0$, so we may take $L = L_{\varepsilon'}$, where $\varepsilon' \in (0, \varepsilon_1)$.

As soon as we know that the set $L$ of possible limits of $\frac{r\partial_r f}{f}$ on $W^\varepsilon$ at $0$ is finite the second part of the proposition follows easily from the standard Łojasiewicz inequality [Lo1].

To show (4.5), it suffices to check that $\frac{|f|}{r^l}$ is bounded from $0$ at $\infty$ on each real analytic curve in $W^\varepsilon$. This follows easily from (4.3). The proof of Proposition 4.2 is now complete. □

We shall need later on the following well-known result.

LEMMA 4.3. *There is a constant $c_f > 0$ such that in a neighborhood of the origin*

(4.6) $$r|\nabla f| \geq c_f |f|.$$

*Proof.* (4.6) is well-known as the Bochnak-Łojasiewicz inequality; see [BL]. It results immediately from the curve selection lemma since $\frac{f}{r|\nabla f|}$ is bounded on each real analytic curve, as again follows easily from (4.3). □

*Remark* 4.4. It seems that the set of characteristic exponents $L \subset \mathbb{Q}^+$ given by Proposition 4.2 is an important invariant of the singularity of $f$. Recall that the Łojasiewicz exponent $\rho_f$ of $f$ is the smallest $\rho$ satisfying (1.1). By Lemma 4.1 and Proposition 4.2, $l \leq (1 - \rho_f)^{-1}$ for $l \in L$. It would be interesting to know whether $(1 - \rho_f)^{-1}$ always belongs to $L$, equivalently whether the Łojasiewicz exponent of $f$ equals $\max_{l_i \in L} \frac{l_i - 1}{l_i}$.

The idea of considering the characteristic exponents $L$ as generalizations of the Łojasiewicz exponent will appear, maybe in a more transparent way, in Corollary 6.5 below.

## 5. Asymptotic critical values

Consider an arbitrary subanalytic $C^1$ function $F$ defined on an open subanalytic set $U$ such that $0 \in \overline{U}$. We say that $a \in \mathbb{R}$ is *an asymptotic critical value of $F$ at the origin* if there exists a sequence $x \to 0$, $x \in U$, such that

(a) $|x||\nabla F(x)| \to 0$ ,

(b) $F(x) \to a$ .

Equivalently we can replace (a) above by

(aa) $|\nabla_\theta F(x)| = |x||\nabla' F(x)| \to 0$ ,

where $\nabla_\theta F$ denotes the gradient of $F$ with respect to spherical coordinates. Indeed, $\nabla_\theta F = r \nabla' F$, so that (a) implies (aa). Suppose that $F(x) \to a$, $|\nabla_\theta F(x)| \to 0$. We have to prove $r \partial_r F \to 0$. If not then there exists a curve, $x = \gamma(r)$, such that on $\gamma$, $|\nabla' F| = o(|\partial_r F|)$ and $|\partial_r F| \geq cr^{-1}$, $c > 0$. In particular, by (4.3),

$$\frac{dF}{dr} \geq \frac{1}{2} cr^{-1}$$

on $\gamma$, so that $F(\gamma(r))$ cannot have a finite limit as $r \to 0$.

PROPOSITION 5.1. *The set of asymptotic critical values is finite.*

*Proof.* Let $X = \{(x,t); F(x) - t = 0\}$ be the graph of $F$. Consider $X$ and $T = \{0\} \times \mathbb{R}$ as a pair of strata in $\mathbb{R}^n \times \mathbb{R}$. Then the (w)-condition of Kuo-Verdier at $(0, a) \in T$ reads

$$1 = |\partial/\partial t (F(x) - t)| \leq C|x||\partial/\partial x(F(x) - t)| = C|x||\nabla F|.$$

In particular, $a \in \mathbb{R}$ is an asymptotic critical value if and only if the condition (w) fails at $(0, a)$. The set of such $a$'s is finite by the genericity of (w) condition; see [V] and [LSW]. □

*Remark* 5.2. The terminology *-an asymptotic critical value-* is motivated by the analogous notion for polynomials $P : \mathbb{R}^n \to \mathbb{R}$ or $P : \mathbb{C}^n \to \mathbb{C}$. We say that $a$ is not an asymptotic critical value of $P$, or equivalently that $P$ satisfies Malgrange's condition at $a$, if there is no sequence $x \to \infty$ such that

(a) $r|\nabla P(x)| \to 0$ ,

(b) $P(x) \to a$ .

For polynomials the set of such values ($a$'s) for which Malgrange's condition fails is finite [Pa3], [KOS]. The proofs there can be easily adapted to the local situation and give alternative proofs of Proposition 5.1. For more on asymptotic critical values see [KOS].

One may ask whether we have an analogue of Łojasiewicz' inequality (1.1) for asymptotic critical values; for instance, whether for an asymptotic critical value $a$ there exist an exponent $\rho_a < 1$, and a constant c, such that

(5.1) $$r|\nabla F| \geq c|F - a|^{\rho_a}.$$

This is not the case in general, but it holds if we approach the singularity "sufficiently slowly."

PROPOSITION 5.3. *Let $F$ be as above and let $a \in \mathbb{R}$. Then for any $\eta > 0$ there exist an exponent $\rho_a < 1$ and constants $c, c_a > 0$, such that (5.1) holds on the set*

$$Z = Z_\eta = \{x \in U; \, |\partial_r F| \leq r^\eta |\nabla F|, |F(x) - a| \leq c_a\}.$$

*Moreover, there exist constants $\delta, \delta' > 0$ such that*

$$Z' = Z'_\delta = \{x \in U; \, r^\delta \leq |F(x) - a| \leq c_a\} \subset Z_{\delta'}.$$

*In particular (5.1) holds on $Z'$.*

*Proof.* For simplicity of notation we suppose $a = 0$. Fix $c_0$ so that $\{|t| \leq c_0\}$ does not contain other asymptotic critical values than 0.

By definition of $Z$

(5.2) $$\frac{\langle \nabla F(x), x \rangle}{|\nabla F(x)||x|} = \frac{\partial_r F}{|\nabla F(x)|} \to 0, \quad \text{as } Z \ni x \to 0.$$

First we show that

(5.3) $$\frac{F(x)}{|\nabla F(x)||x|} \to 0, \quad \text{as } Z \ni x \to 0 \text{ and } F(x) \to 0.$$

It is sufficient to show this on any real analytic curve $\gamma(t)$, such that $\gamma(t) \to 0$ and $F(\gamma(t)) \to 0$ as $t \to 0$. The case $F(\gamma(t)) \equiv 0$ is obvious so we may suppose $F(\gamma(t)) \not\equiv 0$. Note that $\frac{d\gamma/dt}{|d\gamma/dt|}$ and $\frac{\gamma}{|\gamma|}$ have the same limit, so that (5.2) implies

$$\frac{\langle \nabla F(\gamma(t)), d\gamma/dt \rangle}{|\nabla F(\gamma(t))||d\gamma/dt|} \to 0,$$

as $t \to 0$. Hence, $dF/dt = \langle \nabla F, d\gamma/dt \rangle = o(|\nabla F||d\gamma/dt|)$, which gives finally

$$F(x) = o(|\nabla F||x|)$$

along $\gamma$ as required. This demonstrates (5.3) which implies

(5.4) $$\frac{F(x)}{|\nabla F(x)||x|} \to 0, \quad \text{as } x \in Z \text{ and } F(x) \to 0.$$

Indeed, this again has to be checked on any real analytic curve $\gamma(t) \to x_0 \in F^{-1}(0)$, $\gamma(t) \in Z$ for $t \neq 0$. For $x_0 = 0$, it was checked already in (5.3). For $x_0 \neq 0$, $|\gamma|$ has a nonzero limit and (5.4) follows easily by Lemma 4.3 for $F$ at $x_0$.

Finally, (5.4) reads, $\frac{F(x)}{|\nabla F||x|} \to 0$ on $Z$, if $F(x) \to 0$. Hence, by the standard Łojasiewicz inequality [Lo1],

$$\frac{|F(x)|}{|\nabla F(x)||x|} \leq |F|^\alpha, \quad \text{on } Z,$$

for $\alpha > 0$ and sufficiently small. This ends the proof of the first part of Proposition 5.3.

To show the second part we use again the construction from the proof of Proposition 5.1 and the genericity of the Whitney condition (b) for the pair strata $X$ and $T$. Since the Whitney condition (b) is a consequence of the Kuo-Verdier condition (w) for a subanalytic stratification, [V], there is no need to substratify. In particular, for $a'$ not an asymptotic critical value, the Whitney condition (b) implies

$$(5.5) \qquad \frac{\partial_r F}{|\nabla F(x)|} = \frac{\langle \nabla F(x), x \rangle}{|\nabla F(x)||x|} \to 0, \quad \text{as } (x, F(x)) \to (0, a').$$

Let $D = \{t \in T | |t| \leq c_0\}$ so that $D^* = D \setminus \{0\}$ does not contain asymptotic critical values. Then, there is a subanalytic neighborhood $\mathcal{V}$ of $\{0\} \times D^*$ in $\mathbb{R}^n \times \mathbb{R}$, $D^* = D \setminus \{0\}$, such that (5.5) holds for $\mathcal{V} \cap X \ni (x, t) \to T$. Of course, $\mathcal{V}$ can be chosen of the form $\mathcal{V} = \mathcal{V}_\delta = \{(x, t); |x|^\delta \leq |t|, t \in D^*\}$, $\delta > 0$. Now, we may take $Z' = \{x; (x, F(x)) \in \mathcal{V}\}$. Then (5.2) holds as well for $Z' \ni x \to 0$ which implies the existence of $\delta' > 0$ such that $Z' \subset Z_{\delta'}$. □

Next we consider $F$ of the form $F = \frac{f}{r^l}$, $l > 0$, and $U$ the complement of the origin.

PROPOSITION 5.4. *The real number $a \neq 0$ is an asymptotic critical value of $F = \frac{f}{r^l}$ if and only if there exists a sequence $x \to 0$, $x \neq 0$, such that*

(a') $\frac{|\nabla' f(x)|}{|\partial_r f(x)|} \to 0,$

(b) $F(x) \to a.$

*Proof.* Let $x \to 0$ be a sequence in the set $\{x; |\partial_r f| < \varepsilon |\nabla' f|\}$, $\varepsilon > 0$, and such that $F(x) \to a \neq 0$. Then by Lemma 4.3

$$r|\nabla' F| = r\frac{|\nabla' f|}{r^l} \geq c\frac{|f|}{r^l} \geq \frac{1}{2}c|a| > 0.$$

In particular, for such a sequence neither $|x||\nabla F(x)| \to 0$ nor (a') is satisfied.

Thus we may suppose that the sequence $x \to 0$, $F(x) \to a \neq 0$, is in $W^\varepsilon = \{x; |\partial_r f| \geq \varepsilon |\nabla' f|\}$. Then, by Proposition 4.2, we may suppose that

$\frac{|f(x)|}{r^{l_i}}$ is bounded from zero and infinity for an exponent $l_i \in L$. Consequently $l = l_i$. Furthermore, by Proposition 4.2, $|\partial_r f| \sim r^{l-1}$ and

$$\partial_r F = \frac{\partial_r f}{r^l}\left(1 - \frac{lf}{r\partial_r f}\right) = o(r^{-1}), \quad \nabla' F = \frac{\nabla' f}{r^l} = O(r^{-1}).$$

Consequently $a$ is an asymptotic critical value of $F$ if there is a sequence $x \to 0$, $F(x) \to a$, on which

(5.6) $$r|\nabla' F(x)| = \frac{|\nabla' f(x)|}{r^{l-1}} \to 0.$$

Since $|\partial_r f| \sim r^{l-1}$ on $W_l^\varepsilon$, (5.6) is equivalent to (a'). This ends the proof. $\square$

## 6. Estimates on a trajectory

Let $x(s)$ be a trajectory of $\frac{\nabla f}{|\nabla f|}$ defined for $0 \le s < s_0$, $x(s) \to 0$ as $s \to s_0$. In particular, $f(x(s))$ is negative for $s < s_0$. Let $L = \{l_1, \ldots, l_k\}$ denote the set of characteristic exponents of $f$ defined in Proposition 4.2.

Fix $l > 0$, not necessarily in $L$, and consider $F = \frac{f}{r^l}$. Then

(6.1)
$$\begin{aligned}
\frac{dF(x(s))}{ds} &= \langle \frac{\nabla f}{|\nabla f|}, \frac{\nabla' f}{r^l} + \left(\frac{\partial_r f}{r^l} - \frac{lf}{r^{l+1}}\right)\partial_r \rangle \\
&= \frac{1}{|\nabla f| r^l}\left(|\nabla' f|^2 + |\partial_r f|^2 - \frac{lf}{r}\partial_r f\right) \\
&= \frac{1}{|\nabla f| r^l}\left(|\nabla' f|^2 + |\partial_r f|^2\left(1 - \frac{lf}{r\partial_r f}\right)\right).
\end{aligned}$$

PROPOSITION 6.1. *For each $l > 0$ there exist $\varepsilon, \omega > 0$, such that for any trajectory $x(s)$, $F(x(s)) = \frac{f}{r^l}(x(s))$ is strictly increasing in the complement of*

$$\bigcup_{l_i \in L, l_i < l} W_{l_i}^\varepsilon, \quad \text{if} \quad l \notin L,$$

*or in the complement of*

$$W_{-\omega, l} \cup \bigcup_{l_i \in L, l_i < l} W_{l_i}^\varepsilon, \quad \text{if} \quad l \in L,$$

*where in the last case $W_{-\omega, l_i} = \{x \in W_{l_i}^\varepsilon;\ r^{-\omega}|\nabla' f| \le |\partial_r f|\}$.*

*Proof.* If $F$ is not increasing then by (6.1)

$$r|\nabla f|^2 \le lf\partial_r f.$$

Consequently, by (4.6),

$$lf\partial_r f \ge r|\nabla f|^2 \ge c_f|f||\nabla f|,$$

where $c_f$ is the constant of the Bochnak-Łojasiewicz inequality (4.6). In particular, $f\partial_r f$ is positive and

(6.2) $$|\partial_r f| \geq (c_f/l)|\nabla f|.$$

Hence if $F$ is not increasing we are in $W^\varepsilon$ for $\varepsilon = c_f/l$. Recall that $W^\varepsilon = \bigcup W^\varepsilon_{l_i}$ and $\frac{r\partial_r f}{f} \to l_i$ on $W^\varepsilon_{l_i}$. Thus we have three different cases:

- $l < l_i$. Then
$$\left(1 - \frac{lf}{r\partial_r f}\right) \to (1 - l/l_i) > 0.$$
That is to say, $F(x(s))$ is actually increasing in this case.

- $l = l_i$. Then
$$\left(1 - \frac{lf}{r\partial_r f}\right) \to 0$$
on $W^\varepsilon_{l_i}$ and hence is bounded by $\frac{1}{2}r^{2\omega}$, for a constant $\omega > 0$. This means that if $\frac{dF}{ds}$ is negative then $|\partial_r f| \geq r^{-\omega}|\nabla' f|$ as claimed.

- $l > l_i$. Then $F(x(s))$ can be decreasing in $W^\varepsilon_{l_i}$.

This completes the proof of Proposition 6.1. □

PROPOSITION 6.2. *There exist a unique $l = l_i \in L$ and constants: $\varepsilon > 0$ and $0 < c < C < \infty$, such that $x(s)$ passes through $W^\varepsilon_l$ in any neighborhood of the origin and*
$$x(s) \in U_l = \{x;\ c < \frac{|f(x)|}{r^l} < C\}$$
*for $s$ close to $s_0$.*

*Proof.* First we show that the trajectory $x(s)$ passes through $W^\varepsilon$ in any neighborhood of the origin, provided $\varepsilon > 0$ is sufficiently small. Actually any $\varepsilon < c_f(1 - \rho_f)$ would do. Suppose this is not the case. Then, by the proof of Proposition 6.1, $F = \frac{f}{r^l}$ is increasing on the trajectory, for any $l > c_f/\varepsilon > (1 - \rho_f)^{-1}$. Taking into account that $f(x(s))$ is negative we have

(6.3) $$|f(x(s))| \leq C_l r^l,$$

for a $C_l > 0$ which may depend on $l$.

But (6.3) is not possible for $l > (1 - \rho_f)^{-1}$. Indeed, by Łojasiewicz's argument the length of the trajectory between $x(s)$ and the origin is bounded by
$$|s - s_0| \leq c_1 |f(x(s))|^{1-\rho_f}.$$
In particular (6.3) would imply
$$|s - s_0| \leq c_1 C_l r^{l(1-\rho_f)},$$

PROOF OF THE GRADIENT CONJECTURE        779PROOF OF THE GRADIENT CONJECTURE        779

which is not possible for the arc-length parameter $s$ if $l(1 - \rho_f) > 1$. Hence the trajectory $x(s)$ passes through $W^\varepsilon$ in any neighborhood of the origin.

By Proposition 4.2, $W^\varepsilon$ is the finite union of sets $W^\varepsilon_{l_i}$, $l_i \in L$, and each of the $W^\varepsilon_{l_i}$ is contained in a set of the form $U_{l_i} = \{x;\ c < \frac{|f|}{r^{l_i}} < C\}$. We fix $c$ and $C$ common for all $l_i$. The $U_{l_i}$'s are mutually disjoint as germs at the origin.

Fix one of the $l_i$'s and consider $F = \frac{f}{r^{l_i}}$. By Proposition 6.1, $F(x(s))$ is strictly increasing on the boundary of $U_{l_i}$, that is to say if

$$x(s) \in \partial^+ U_i = \{x;\ f(x) = -Cr^{l_i}\},$$

then the trajectory enters $U_{l_i}$ and if

$$x(s) \in \partial^- U_{l_i} = \{x;\ f(x) = -cr^{l_i}\},$$

then the trajectory leaves $U_{l_i}$, of course definitely. This ends the proof. $\square$

*Remark* 6.3. Note that the constants $\varepsilon$ and $c, C$ of Proposition 6.2 can be chosen independent of the trajectory. Indeed, by the proof of Proposition 6.1, we may choose, for instance, $\varepsilon \leq \frac{1}{2} c_f (1 - \rho_f)$. Then, by Remark 4.4, $\varepsilon \leq \frac{1}{2} c_f / l$ for any $l \in L$. Now the constants $c$ and $C$ are given by (4.5) and we fix as $\omega > 0$ any exponent which satisfies

$$(6.4) \qquad \left|1 - \frac{lf}{r \partial_r f}\right| \leq \frac{1}{2} r^{2\omega}$$

on $W^\varepsilon_l$ for each $l \in L$.

Let us list below some of the bounds satisfied on $W^\varepsilon_l$ and $U_l$. Recall that, by construction, $U_l \supset W^\varepsilon_l$. If $\varepsilon \leq \frac{1}{2}(c_f/l)$, as we have assumed, then by (6.2), we have, away from $W^\varepsilon_l$,

$$(6.5) \qquad r|\nabla f|^2 \geq 2 l f \partial_r f,$$

which gives

$$\frac{dF(x(s))}{ds} \geq \frac{|\nabla f|}{2r^l} \geq \frac{c_f |f|}{2 r^{l+1}}.$$

Hence on $U_l \setminus W^\varepsilon_l$

$$(6.6) \qquad \frac{dF(x(s))}{ds} \geq \frac{c_f |f|}{2 r^{l+1}} \geq c' r^{-1}$$

for a universal constant $c' > 0$. Also by Section 4 we have easily

$$|\nabla f| \geq c_1 r^{l-1} \quad \text{on } U_l,\ \text{for } c_1 > 0,$$

$$|\nabla f| \sim \partial_r f \sim r^{l-1} \quad \text{on } W^\varepsilon_l.$$

From now on we shall assume $\varepsilon$ and $\omega$ fixed.

We shall show in the proposition below that $F = \frac{f(x(s))}{r^l}$ has a limit as $s \to s_0$. For this we use an auxiliary function $F - r^\alpha$.

PROPOSITION 6.4. *For $\alpha < 2\omega$, the function $g = F - r^\alpha$ is strictly increasing on the trajectory $x(s)$. In particular $F(x(s))$ has a nonzero limit*

$$F(x(s)) \to a_0 < 0, \quad as\ s \to s_0.$$

*Furthermore, $a_0$ must be an asymptotic critical value of $F$ at the origin.*

*Proof.* By Proposition 6.1, $F$ is increasing on $|\partial_r f| \leq r^{-\omega}|\nabla' f|$. On the other hand

$$\frac{d(-r^\alpha)}{ds} = -\alpha r^{\alpha-1} \frac{\partial_r f}{|\nabla f|},$$

and hence $-r^\alpha$ is increasing if $\partial_r f$ is negative which is the case on $W_l^\varepsilon$ (on $W_l^\varepsilon$, $f \partial_r f > 0$ and $f < 0$ on the trajectory). We consider three different cases:

*Case 1.* $r^{-\omega}|\nabla' f| \leq |\partial_r f|$. That is, we are in $W_{-\omega,l}$ of Proposition 6.2.

Then, in particular, we are in $W_l^\varepsilon$ and (6.4) holds. Moreover $|\partial_r f| \sim |\nabla f| \sim r^{l-1}$ (see Remark 6.3) and $\partial_r f$ is negative. Consequently

$$\left|\frac{dF(x(s))}{ds}\right| \leq \frac{1}{|\nabla f| r^l}\left(|\nabla' f|^2 + |\partial_r f|^2 r^{2\omega}\right) \leq C_1(r^{2\omega-1}),$$

$$\frac{d(-r^\alpha)}{ds} = -\alpha r^{\alpha-1}\frac{\partial_r f}{|\nabla f|} \geq (\alpha/2) r^{\alpha-1}.$$

Thus $\frac{d(-r^\alpha)}{ds}$ is dominant and $g$ is increasing on the trajectory.

*Case 2.* $\varepsilon|\nabla' f| \leq |\partial_r f| < r^{-\omega}|\nabla' f|$. That is, we are in $W_l^\varepsilon \setminus W_{-\omega,l}$.
Then both $F$ and $-r^\alpha$ are increasing.

*Case 3.* $|\partial_r f| < \varepsilon|\nabla' f|$. That is, we are in $U_l \setminus W_l^\varepsilon$.
By Remark 6.3,

$$\frac{dF(x(s))}{ds} \geq c' r^{-1}.$$

On the other hand

$$\left|\frac{dr^\alpha}{ds}\right| = \left|\alpha r^{\alpha-1}\frac{\partial_r f}{|\nabla f|}\right| \leq r^{\alpha-1},$$

so that $g$ is increasing.

Finally, since $g(x(s))$ is increasing, negative and bounded from zero on $U_l$, it has the limit $a_0 < 0$. We shall show that $a_0$ is an asymptotic critical value of $F$.

Suppose that, contrary to our claim, $F(x(s)) \to a_0$ and $a_0$ is not an asymptotic critical value of $F$ at the origin. Then, by Proposition 5.4, there is $\tilde{c} > 0$ such that

$$|\nabla' f(x(s))| \geq \tilde{c}|\partial_r f(x(s))|,$$

for $s$ close to $s_0$. Hence on $W_l^\varepsilon$

(6.7) $$\frac{dF}{ds} = \frac{|\nabla' f|^2}{r^l |\nabla f|} + \frac{|\partial_r f|^2}{r^l |\nabla f|}\left(1 - \frac{lf}{r \partial_r f}\right) \geq c' \frac{1}{r}.$$

A similar bound holds on $U_l \setminus W_l^\varepsilon$ by (6.6).

But (6.7) is not possible since $r \leq s_0 - s$. Indeed, the inequality (6.7) implies $\frac{dF}{ds} \geq c \frac{1}{s_0 - s}$ with the right-hand side not integrable which contradicts the fact that $F$ is bounded on the trajectory (Proposition 6.1). This ends the proof. □

COROLLARY 6.5. *Let $\sigma(s)$ denote the length of the trajectory between $x(s)$ and the origin. Then*

$$\frac{\sigma(s)}{|x(s)|} \to 1 \text{ as } s \to s_0.$$

*Proof.* We show that if $x(s)$ is entirely contained in $U_{l,\tilde{c},\tilde{C}} = \{x | 0 < \tilde{c} < \frac{|f|}{r^l} < \tilde{C} < \infty\}$ then

(6.8) $$\sigma(s) \leq [\tilde{C}/\tilde{c}]^{1/l} r + o(r),$$

where $r$ denotes $|x(s)|$. Then the corollary follows directly from the existence of the limit of $F(x(s))$ as $s \to s_0$.

In order to establish (6.8) we follow in detail the computation of the Łojasiewicz argument. First we note that on $U_{l,\tilde{c},\tilde{C}}$,

$$r|\nabla f| \geq l|f| - o(r^l),$$

which can be checked on any real analytic curve as in (4.3). Then

$$\frac{df}{ds} = |\nabla f| \geq \frac{l|f|}{r} - o(r^{l-1}) \geq l\tilde{c}^{1/l} |f|^{\rho_l} - o(r^{l-1}),$$

where $\rho_l = \frac{l-1}{l}$. Hence, by the Łojasiewicz argument

$$\sigma(s) \leq [l\tilde{c}^{1/l}(1 - \rho_l)]^{-1} |f|^{1-\rho_l} \leq [\tilde{C}/\tilde{c}]^{1/l} r + o(r),$$

as required. □

## 7. Proof of main theorem

THEOREM 7.1. *Let $x(s)$ be a trajectory of $\frac{\nabla f}{|\nabla f|}$, $x(s) \to 0$ as $s \to s_0$. Denote by $\tilde{x}(s)$ the radial projection of $x(s)$ onto the unit sphere, $\tilde{x}(s) = \frac{x(s)}{|x(s)|}$. Then $\tilde{x}(s)$ is of finite length.*

*Proof.* By Section 6 there exists a unique $l \in L$, $L$ given by Proposition 4.2, and a unique asymptotic critical value $a < 0$ of $F = \frac{f}{r^l}$, such that

$$F(x(s)) \to a, \quad \text{as } s \to s_0.$$

Thus we may suppose that $x(s)$ is entirely included in

$$U_{l,a} = \{x \in U_l; |F(x) - a| \leq c_a\},$$

where $c_a$ is as given by Proposition 5.3.

We use the arc-length parametrization $\tilde{s}$ of $\tilde{x}(s)$ given by

$$\frac{ds}{d\tilde{s}} = \frac{r|\nabla f|}{|\nabla' f|}.$$

Reparametrize $x(s)$ using $\tilde{s}$ as parameter. Then

$$\frac{dF}{d\tilde{s}} = \frac{r|\nabla f|}{|\nabla' f|} \frac{dF}{ds} = \frac{1}{|\nabla' f| r^{l-1}} \left( |\nabla' f|^2 + |\partial_r f|^2 \left(1 - \frac{lf}{r \partial_r f}\right) \right),$$

$$\frac{d(-r^\alpha)}{d\tilde{s}} = -\alpha r^\alpha \frac{\partial_r f}{|\nabla' f|}.$$

Let $\omega$ be the exponent given by (6.4) and let $\rho_1 = \rho_a$ be the exponent given by Proposition 5.3 for $\eta < \omega$.

LEMMA 7.2.

(7.1) $$\frac{d(F(x(s)) - a)}{d\tilde{s}} \geq c|F(x(s)) - a|^{\rho_1}$$

on $\{x; |\partial_r f| \leq r^{-\eta} |\nabla' f|\}$ and for $c > 0$.

*Proof.* If $|\partial_r f| < \varepsilon |\nabla' f|$ (Case 3 of the proof of Proposition 6.4), then

(7.2) $$\frac{dF}{d\tilde{s}} \geq c' r^{-1} \frac{r|\nabla f|}{|\nabla' f|} \geq c' > 0,$$

and (7.1) follows trivially.

So suppose $\varepsilon |\nabla' f| \leq |\partial_r f| \leq r^{-\eta} |\nabla' f|$. Recall that after Remark 6.3, $|\nabla f| \sim |\partial_r f| \sim r^{l-1}$. Hence, by $\eta < \omega$ and (6.4),

$$r|\nabla' F| = \frac{|\nabla' f|}{r^{l-1}} \geq c_1 r^\eta \gg r^{2\omega - \eta} \sim r^{-\eta} |\frac{\partial_r f}{r^{l-1}} \left(1 - \frac{lf}{r \partial_r f}\right)| = r^{1-\eta} |\partial_r F|.$$

Consequently,

$$\frac{dF}{d\tilde{s}} \geq \frac{|\nabla' f|}{r^{l-1}} + \frac{|\partial_r f|^2}{r^{l-1} |\nabla' f|} \left(1 - \frac{lf}{r \partial_r f}\right) \geq \frac{1}{2} r |\nabla F|$$

and the statement follows from Proposition 5.3. □

Thus, in the complement of $W_{-\eta,l} = \{x \in W_l^\varepsilon \mid r^{-\eta}|\nabla' f| \leq |\partial_r f|\}$, $(F-a)$ can be used as a control function in the sense of Łojasiewicz. To get a control function which works everywhere on $U_l$ we consider

$$(7.3) \qquad g = (F-a) - r^\alpha$$

for any $0 < \alpha < \omega$. Our aim is to show the existence of $\xi$, $\xi < 1$, such that

$$(7.4) \qquad \frac{dg}{d\tilde{s}} \geq c|g|^\xi \quad \text{for } c > 0$$

on $U_{l,a}$. Note that $-r^\alpha$ is itself a good control function on $W_{-\delta,l} = \{x \in W_l^\varepsilon; r^{-\delta}|\nabla' f| \leq |\partial_r f|\}$, where we may chose $\delta > 0$ arbitrarily. Indeed, on $W_{-\delta,l}$

$$\frac{d(-r^\alpha)}{d\tilde{s}} = -\alpha r^\alpha \frac{\partial_r f}{|\nabla' f|} \geq \alpha r^{\alpha-\delta} \geq \alpha (r^\alpha)^\rho,$$

for $\rho \geq 1 - \frac{\delta}{\alpha}$.

First we show that replacing $F$ by $g$ does not destroy the property of being a control function in the complement of $W_{-\eta,l}$.

LEMMA 7.3. *For any $\alpha > 0$, $g = (F-a) - r^\alpha$ is a control function on $U_{l,a} \setminus W_{-\eta,l}$.*

*Proof.* By (7.2) there is no problem on $U_l \setminus W_l^\varepsilon$. Of course, there is no problem in the area where both summands are good control functions. Thus it suffices to consider $W_l^\varepsilon \setminus W_{-\delta,l}$, where we may choose $\delta > 0$ arbitrarily small. But on this set $\frac{dF}{d\tilde{s}} \geq \frac{1}{2}r|\nabla F| \geq cr^\delta$ and hence

$$\frac{dg}{d\tilde{s}} = \frac{d(F-a)}{d\tilde{s}} - \frac{dr^\alpha}{d\tilde{s}} \geq c(r^\delta + r^\alpha) \geq cr^\delta \geq c(r^\alpha)^\rho, \quad c > 0,$$

if we set, for example: $\delta < \alpha$ and $\rho \geq 1 - \frac{\delta}{\alpha}$ which proves the lemma. □

Thus it remains to consider $W_{-\eta,l} = \{x \in W_l^\varepsilon; r^{-\eta}|\nabla' f| \leq |\partial_r f|\}$. On $W_{-\eta,l}$ the radial part of the gradient is much bigger than the spherical one and therefore $-r^\alpha$ can be taken as a control function. Moreover, as we show below, the derivative of $-r^\alpha$, $\alpha \leq \eta$, dominates the derivative of $F-a$ along the trajectory. On $W_{-\eta,l}$

$$\frac{dF}{d\tilde{s}} = \frac{|\nabla' f|}{r^{l-1}} + \frac{|\partial_r f|}{r^{l-1}} \frac{|\partial_r f|}{|\nabla' f|} O(r^{2\omega}) = O(r^\eta) + O(r^{2\omega-\eta})$$

and

$$\frac{d(-r^\alpha)}{d\tilde{s}} = \alpha r^\alpha \frac{|\partial_r f|}{|\nabla' f|} \geq \alpha r^{\alpha-\eta} \geq \text{const} > 0$$

and consequently $\frac{dg}{d\tilde{s}} \geq \text{const} > 0$. Together with Lemma 7.3 this implies that $g$ is a control function everywhere on $U_{l,a}$, in the sense that there exists $\xi < 1$ such that (7.4) holds on $U_{l,a}$.

Now the theorem follows directly by Łojasiewicz's argument, recalled below for completeness. We proved that (7.4) holds on $U_{l,a} = \{x \in U_l; |F(x) - a| \leq c_a\}$. Hence, there is $c_0$ such that

$$\frac{d|g|^{1-\xi}}{d\tilde{s}} \leq c_0 < 0$$

on $U_{l,a}$. Let $s_1 < s_0$ be such that $x(s) \in U_{l,a}$ for $s_1 \leq s < s_0$. Recall that $g(x(s)) \to 0$ as $s \to s_0$ and $g < 0$ on the trajectory. Hence, by integration,

(7.5) $$\text{length}\{\tilde{x}(s); s_1 \leq s \leq s_0\} \leq c_0^{-1} |g(x(s_1))|^{1-\xi}$$

which is finite. This ends the proof. □

## 8. Miscellaneous

In the first part of this section we show that there exists a uniform bound on the length of the spherical part of the trajectory.

In the second part we derive more information on the way the trajectory $x(s)$ of the gradient approaches its limit point $x_0 = \lim_{s \to s_0} x(s)$. Theorem 7.1 gives the existence of a tangent line to $x(s)$ at $x_0$. Denote this line by $P$. One may expect that, for $s$ close to $s_0$, $x(s)$ stays in a cuspidal (horn) neighborhood of $P$, in other words that the distance from $x(s)$ to $P$ can be bounded by a positive power of $r = |x(s) - x_0|$. As shown in an example below, in general this is not the case. On the other hand one may show that this distance is bounded by a continuous function $\alpha(r) \to 0$ which can be expressed in terms of real exponential and logarithmic functions (i.e. $\varphi$ belongs to the o-minimal structure defined by polynomials and the exponential function).

THEOREM 8.1.  *There exist $M > 0$ and $r_0 > 0$ such that for any trajectory $x(s)$, $0 \leq s < s_0$, of $\frac{\nabla f}{|\nabla f|}$, which tends to the origin, and which is contained entirely in the ball $B(0, r_0) = \{x; |x| \leq r_0\}$,*

(8.1) $$\text{length}\{\tilde{x}(s)\} \leq M.$$

*Proof.* It is easy to check that the proof of Theorem 7.1 gives such a universal bound as long as $x(s)$ is contained in $U_{l,a}$. But, in general, the moment when the trajectory falls into $U_{l,a}$ depends on the trajectory. On the other hand, it is not difficult to generalize the arguments of Section 7 to show that the control function $g = F - r^\alpha$, for $\alpha > 0$ sufficiently small, gives a universal bound on the length of the part of $\tilde{x}(s)$ which is contained in $U_l$. We have an even stronger statement.

LEMMA 8.2. *There are $\delta > 0$ and a constant $M_1$ such that for any $l \in L$, if*

$$(8.2) \qquad x(s) \in U_l^\delta = \{x; -r^{-\delta} \leq \frac{f}{r^l} \leq -r^\delta\},$$

*for $s_1 \leq s < s_2$, then* $\text{length}\{\tilde{x}(s); s_1 \leq s < s_2\} \leq M_1$.

*Proof.* We shall always assume that the $U_l^\delta$'s, for $l \in L$, are mutually disjoint, which of course we may do in a small neighborhood of the origin and for $\delta$ sufficiently small. In particular, this implies that that $W_{l_i}^\varepsilon \cap U_{l_j}^\delta = \emptyset$ for $l_i \neq l_j$.

First we show the statement on $U_l = \{x; c < \frac{|f(x)|}{r^l} < C\}$. Note that $F = \frac{f}{r^l}$ satisfies (5.1) on $U_l \setminus \bigcup_{a \in A_l} U_{l,a}$, where $A_l$ denotes the set of all negative asymptotic critical values of $F$. Hence, by the proof of Theorem 7.1, $g = \frac{f}{r^l} - r^\alpha$, $\alpha > 0$ sufficiently small, satisfies

$$(8.3) \qquad \frac{dg}{d\tilde{s}} \geq \prod_{a_i \in A_l} c|g - a_i|^\xi,$$

for $\xi < 1$ and $c > 0$.

In particular, $g(x(s))$ is increasing in $U_l$. Consider

$$h(x) = \sum_{a_i \in A_l} (g(x) - a_i)^{1-\xi'},$$

where $\xi' < 1$ is chosen so that $\xi \leq \xi' < 1$ and $1 - \xi'$ is the quotient of two odd integers. Then each summand of $h$ is increasing on the trajectory and (8.3) implies

$$\frac{dh}{d\tilde{s}} \geq c_0 > 0,$$

which gives by integration the required bound on the length of the part of the trajectory which is in $U_l$:

$$\text{length}\{\tilde{x}(s); x(s) \in U_l\} \leq c_0^{-1}[\max\{h(x); x \in U_l\} - \min\{h(x); x \in U_l\}]$$

(recall that $x(s)$ may cross $U_l$ only once; see the proof of Proposition 6.2).

Secondly, we consider $U_l^{\delta,+} = \{x; -c \leq \frac{f}{r^l} \leq -r^\delta\}$, $\delta > 0$ to be specified below. For this we use Proposition 5.3 for $a = 0$ which implies the existence of $\delta > 0$ and $\rho_0 < 1$ such that

$$r|\nabla F| \geq c_1 |F|^{\rho_0}$$

on $|F| \geq r^\delta$. Note that, since we are away from $W^\varepsilon$, (6.5) holds, and hence

$$\frac{dF}{d\tilde{s}} = \frac{1}{|\nabla' f|r^{l-1}}\left(|\nabla' f|^2 + |\partial_r f|^2\left(1 - \frac{lf}{r\partial_r f}\right)\right) \geq \frac{1}{2}r|\nabla F| \geq \frac{1}{2}c_2|F|^{\rho_0}.$$

Thus we may use Łojasiewicz's argument. Consequently the length of $\tilde{x}(s)$, corresponding to the part of $x(s)$ contained in $U_l^{\delta,+}$, can be bounded by a universal constant.

Finally, we consider $U_l^{\delta,-} = \{x; -r^{-\delta} \leq \frac{f}{r^l} \leq -C\}$. Now $g = -F^{-1}$ can be used as a control function. Note that $\nabla g = F^{-2}\nabla F$ so in the set $\{x; r^\delta \leq g \leq C^{-1}\}$, $g$ has no asymptotic critical values except maybe zero. We use again Proposition 5.3 for $g$ and $a = 0$

$$(8.4) \qquad \frac{dg}{d\tilde{s}} = g^{-2}\frac{dF}{d\tilde{s}} \geq \frac{1}{2}rg^{-2}|\nabla F| \geq \frac{1}{2}r|\nabla g| \geq c_3 g^{\rho_\infty},$$

on the set $|g| \geq r^\delta$ and for $\rho_\infty < 1$, $c_3 > 0$. Again by Łojasiewicz's argument we get a required bound on $U_l^{\delta,-}$. This ends the proof of the lemma. □

Note that the trajectory $x(s)$ may cross each of the sets $U_l^{\delta,-}$, $U_l$, and $U_l^{\delta,+}$ only once. Indeed, for $U_l$ we proved it previously, in Proposition 6.2. A similar argument works for $U_l^{\delta,-}$ and $U_l^{\delta,+}$. For instance, by Proposition 6.1, $\frac{f}{r^{l-\delta}}$ is strictly increasing on the boundary of $U_l^{\delta,-}$. Hence the trajectory may enter $U_l^{\delta,-}$ only once.

Now we consider the complement $U^0$ of $\bigcup_{l \in L} U_l^\delta$. On this set the length of $\tilde{x}(s)$ is easy to bound since, as we prove below, the spherical component of the movement of $x(s)$ dominates the radial one and we may use $f$ as a control function. More precisely we have the following lemmas.

LEMMA 8.3. *There exists an exponent $\delta_1 > 0$ such that in the complement $U^0$ of $\bigcup_{l \in L} U_l^\delta$*

$$|\nabla' f| \geq r^{-\delta_1}|\partial_r f|.$$

*Proof.* Let $\gamma(t)$ be an analytic curve in $U^0 \cup \{0\}$, such that $\gamma(t) \to 0$ as $t \to 0$. Then

$$\frac{\partial_r f}{|\nabla' f|}(\gamma(t)) \to 0 \quad \text{as } t \to 0.$$

Indeed, if this is not the case then $\gamma(t) \in W^{\varepsilon'}$, for $\varepsilon'$ sufficiently small. But, by Proposition 4.2, $W^{\varepsilon'}$ is contained in $\bigcup_{l \in L} U_l^\delta$ at least in a neighborhood of the origin, which contradicts the assumptions. Hence $\frac{\partial_r f}{|\nabla' f|}(x) \to 0$ as $U^0 \ni x \to 0$. This implies, by the standard Łojasiewicz inequality, the existence of $\delta_1 > 0$ such that

$$\frac{\partial_r f}{|\nabla' f|} \leq r^{\delta_1},$$

as required. □

The following lemma, which we will prove directly, is a version of Proposition 5.3 applied to $f$ and zero as the (only) critical value of $f$.

LEMMA 8.4. *For each $\delta_1 > 0$ there exists $\tilde{\rho} < 1$ such that on the set $V_{\delta_1} = \{x; |\nabla' f| \geq r^{-\delta_1}|\partial_r f|\}$,*

$$r|\nabla f| \geq |f|^{\tilde{\rho}}.$$

*Proof.* The lemma follows again by direct application of the standard Łojasiewicz inequality and the curve selection lemma. The quotient $\frac{f}{r|\nabla f|}(x) \to 0$ as $V_{\delta_1} \ni x \to x_0 \in f^{-1}(0)$, as can be checked on any real analytic curve. Hence this quotient is bounded by a power of $f$, as required. □

Now we are in a position to complete the proof of Theorem 8.1. By Lemmas 8.3 and 8.4, we have in $U^0$

$$\frac{df}{d\tilde{s}} = \frac{r|\nabla f|^2}{|\nabla' f|} \geq r|\nabla f| \geq |f|^{\tilde{\rho}}.$$

So each segment of $x(s)$, $s \in [s_1, s_2]$, included in $U^0$, has the length of its radial projection bounded by $(1-\tilde{\rho})^{-1}[|f(x(s_1))|^{1-\tilde{\rho}} - |f(x(s_2))|^{1-\tilde{\rho}}]$. Hence the total length of the radial projection of the part of $x(s)$ which is in $U^0$ is bounded by $(1-\tilde{\rho})^{-1}|f(x(0))|^{1-\tilde{\rho}}$.

This ends the proof. □

Note that we cannot require that the constant $M$ of the theorem above tends to 0 as $r_0 \to 0$. Indeed, consider $f : \mathbb{R}^2 \to \mathbb{R}$ given by $f(x,y) = -\frac{1}{2}(x^2 + ay^2)$, $a > 0, a \neq 1$. Then the best possible uniform bound on the length of $\tilde{x}(s)$, even for $r_0$ arbitrarily small, is $\pi/2$.

*Example* 8.5. Consider the trajectories of

$$f(x,y) = -\frac{1}{4}(x^2+y^2)^2 - \frac{1}{3}xy^3.$$

In polar coordinates $f(r,\theta) = -\frac{1}{4}r^4(1 + \frac{4}{3}\cos\theta\sin^3\theta)$,

$$\frac{\partial f}{\partial r} = -r^3(1 + \frac{4}{3}\cos\theta\sin^3\theta), \qquad r^{-2}\frac{\partial f}{\partial \theta} = r^2(\frac{1}{3}\sin^4\theta - \cos^2\theta\sin^2\theta).$$

Therefore, by (2.2), the trajectories of the gradient in polar coordinates satisfy

$$\frac{d\theta}{dr} = \frac{\theta^2(1 + \sigma(\theta))}{r},$$

where $\sigma$ is an analytic function at 0 and $\sigma(\theta) = O(\theta^2)$. By integration we get a solution $\theta(r)$ of the form

$$\theta(r) \simeq (-\ln r)^{-1}.$$

In particular $\theta(r) \to 0$ slower than any positive power of $r$.

Note also that in this example the control function of the proof of the main theorem $F - a$, which as we proved tends to 0 along the trajectory, goes to 0 more slowly than any positive power of $r$.

PROPOSITION 8.6. *Under the hypotheses of Theorem* 7.1 *there exists a continuous function* $\varphi(r) : (\mathbb{R}, 0) \to (\mathbb{R}, 0)$, (*which depends on the trajectory*), *of the form* $\varphi(r) = (-\ln r)^{-\delta}$, $\delta > 0$, *such that*

(8.5) $$\text{length}\{\tilde{x}(s), s_1 \leq s < s_0\} \leq \varphi(|x(s_1)|).$$

*Proof.* We follow the notation of previous sections. Let $g = (F - a) - r^\alpha$ be the control function of the proof of Theorem 7.1.

LEMMA 8.7. *There is an exponent* $\mu$ *such that*

(8.6) $$r \frac{dg}{ds} \geq |g|^\mu$$

*on the trajectory.*

*Proof.* We follow the cases of the proof of Proposition 6.4.

In Case 3, that is, in the complement of $W_l^\varepsilon$,

$$r \frac{dg}{ds} \geq c \frac{dg}{d\tilde{s}} \geq \text{const} > 0;$$

see (7.2), and the lemma follows.

On $W_l^\varepsilon$, $\partial_r f$ is negative and consequently

$$r \frac{d(-r^\alpha)}{ds} = -\alpha r^\alpha \frac{\partial_r f}{|\nabla f|} \geq c r^\alpha,$$

where $c > 0$. That is (8.6) holds for $-r^\alpha$. Therefore, in both Cases 1 and 2 of the proof of Proposition 6.4,

$$r \frac{dg}{ds} \geq -(\alpha/2) r^\alpha \frac{\partial_r f}{|\nabla f|} \geq c' r^\alpha$$

with $c' > 0$. Therefore the lemma certainly holds on the sets defined by $(F - a) \leq r^\delta$, $\delta > 0$.

We suppose that $(F - a) \geq r^\delta$ with $\delta > 0$ sufficiently small. We observe that in Proposition 5.3, for each $\delta' > 0$ there is $\delta > 0$ such that $Z'_\delta \subset Z_{\delta'}$ (we suppose $c_a$ of the proposition small). Indeed, it follows from the fact that the sets $\mathcal{V}_\delta$, $\delta > 0$ (see the last paragraph of the proof of Proposition 5.3) form a base of subanalytic neighborhoods of $D^*$. Hence for both $\delta > 0$ and $\delta' > 0$ sufficiently small, by (6.1) and (5.1),

$$r \frac{dF}{ds} \geq \frac{c|\nabla' f|^2}{r^{2(l-1)}} \geq c_1 \Big(\frac{dF}{d\tilde{s}}\Big)^2 \geq c_2 |F - a|^{2\rho_a},$$

$c, c_1, c_2 > 0$. Thus both summands of $g$, $F - a$ and $-r^\alpha$ have positive derivatives and for both of them the inequality of type (8.6) holds. Hence it holds for $g$. This ends the proof of Lemma 8.7. □

We shall use as parameter $\tau = s_0 - s$. Also we replace $g$ by its absolute value $h$. Now the Lemma 8.7 reads: there is an exponent $\mu$ such that

$$\tau \frac{dh}{d\tau} \geq h^\mu.$$

Without loss of generality we may suppose $\mu > 1$. By integration we get on the trajectory

(8.7) $$h^{\mu-1} \leq (C - (\mu-1)\ln\tau)^{-1}.$$

A similar bound holds if we replace $\tau$ by $r$, $r \simeq \tau$ by Corollary 6.5. Now the proposition follow easily from (7.5). □

## 9. The Riemannian case

Let $(M, g)$ will be an analytic manifold equipped with an analytic riemannian metric $g$. Letting $f : M \to \mathbb{R}$ be an analytic function, we denote by $\nabla_g f$ the gradient of $f$ with respect to $g$. Let $x : [0, s_0) \to M$ be a trajectory of $\nabla_g f$, i.e., a maximal solution of

$$x' = \frac{\nabla_g f}{|\nabla_g f|}, \quad x(0) = p_0$$

for some $p_0 \in M$. Suppose that $s \to s_0$; then either $x(s)$ escapes from $M$ or has a limit point $x_0$ in $M$ (this means that the trajectory passes by any neighborhood of $x_0$). Suppose that we are in the second case and assume that $f(x_0) = 0$.

It is easily seen that the Łojasiewicz inequality $|\nabla_g f| \geq c|f|^\rho$, $\rho < 1$, holds in a neighborhood of $x_0$. Thus we deduce (as in the introduction) that: $x_0 = \lim_{s \to s_0} x(s)$; moreover $s_0 < \infty$ and $x_0$ is a critical point of $f$. Now it makes sense to ask if the curve $x(s)$ has a tangent at the limit point $x_0$; this is exactly the formulation of the gradient conjecture in the Riemannian case. We can put it also in the following way: let $\varphi : U \to \mathbb{R}^n$, $\varphi(x_0) = 0$ be a chart around point $x_0$. Does the limit of secants

$$\lim_{s \to s_0} \frac{\varphi(x(s))}{|\varphi(x(s))|}$$

exist?

Here $||$ stands for the Euclidean norm, but of course it can be any norm on $\mathbb{R}^n$. By the definition of a manifold (compatibility of charts) the answer does not depend on the choice of the chart $\varphi$. Equivalently we can state it as in [Mu], assuming that $M = \mathbb{R}^n$ is equipped with some Riemannian metric $g$ and $\varphi$ is simply the identity on $\mathbb{R}^n$. So the gradient conjecture asks if $\lim_{s \to s_0} \frac{x(s)}{|x(s)|}$ exists, where $x(s)$ is a trajectory of $\nabla_g f$, $x(s) \to 0$ as $s \to s_0$. We prove below that actually:

THEOREM 9.1.  *The gradient conjecture holds in the Riemannian case.*

To prove it we first embed isometrically and analytically $(M, g)$ into $\mathbb{R}^N$ equipped with the Euclidean metric. Next, we prove that the gradient conjecture (existence of the limit of secants) holds for analytic submanifolds with the metric induced by the Euclidean metric. In fact, by using analytic tubular neighborhoods we reduce the problem to the main Theorem 7.1 in the Euclidean case for an open set in $\mathbb{R}^N$.

We shall need a local version of the isometric embedding theorem in the analytic case due to M. Janet [Ja] and E. Cartan [Ca].

THEOREM 9.2 (Janet-Cartan).   *Let $(M, g)$ be an analytic Riemannian manifold. Then, for any $x_0 \in M$ there exist a neighborhood $U$ and an analytic embedding $\sigma : U \to \mathbb{R}^N$ such that $\sigma(U)$ is a submanifold of $\mathbb{R}^N$ and $\sigma$ is an isometry between Riemannian manifolds $(U, g)$ and $(\sigma(U), e)$, where $e$ is the Riemannian metric on $\sigma(U)$ induced by the Euclidean metric on $\mathbb{R}^N$.*

So we may assume that $M$ is a an analytic submanifold of $\mathbb{R}^N$ and $g$ is induced by the Euclidean metric on $\mathbb{R}^N$.

PROPOSITION 9.3.   *Suppose that $f : M \to \mathbb{R}$ is an analytic function and $x(s) \in M$ is a trajectory of $\frac{\nabla_g f}{|\nabla_g f|}$, $x(s) \to 0$ as $s \to s_0$; then $\tilde{x}(s) = \frac{x(s)}{|x(s)|}$ is of finite length.*

*Proof.* Let $\pi : U \to M$ be an analytic tubular neighborhood of $M$ in $\mathbb{R}^N$ with respect to the Euclidean metric. So in particular for any $x \in M$ the fiber $\pi^{-1}(x)$ can be seen as an open subset of the normal space to $M$ at $x$. Let us put $h = f \circ \pi$. Take $\nabla h(x)$, the gradient of $h$ with respect to the Euclidean metric; then

$$\nabla h(x) = \nabla_g f(x) \quad \text{for any } x \in M.$$

Hence our trajectory $x(s)$ is also a trajectory of $\frac{\nabla h}{|\nabla h|}$, so the proposition follows from Theorem 7.1. This also proves Theorem 9.1; indeed we can take as chart $\varphi$ the orthogonal projection of $M$ on the tangent space to $M$ at 0.   □

JAGIELLONIAN UNIVERSITY, KRAKÓW, POLAND, AND
UNIVERSITÉ DE SAVOIE, LE BOURGET-DU-LAC, FRANCE
*E-mail address*: kurdyka@univ-savoie.fr

WARSAW UNIVERSITY, WARSAW, POLAND
*E-mail address*: tmostows@ghost.mimuw.edu.pl

UNIVERSITÉ D'ANGERS, ANGERS, FRANCE
*E-mail address*: parus@tonton.univ-angers.fr